\documentclass{gtart}

\def\ifplaintex{\expandafter\ifx\csname documentclass\endcsname\relax}

\def\gtp{{\mathsurround=0pt\it $\cal G\mskip-2mu$eometry \&\ 
$\cal T\!\!$opology $\cal P\!$ublications}}  % GT publications

\def\recd{{\small Received:\qua\receiveddate\ifx\reviseddate\relax
\else\qquad Revised:\qua\reviseddate\fi\par}} 

\def\lognumber#1{\def\thelognumber{#1}}
\def\volumenumber#1{\def\thevolumenumber{#1}}
\def\volumeyear#1{\def\thevolumeyear{#1}}
\def\papernumber#1{\def\thepapernumber{#1}}
\def\pagenumbers#1#2{\def\startpage{#1}\def\finishpage{#2}}
\def\published#1{\def\publishdate{#1}}

\def\received#1{\def\receiveddate{#1}}
\def\revised#1{\def\reviseddate{#1}}
\def\accepted#1{\def\accepteddate{#1}}

\def\asciiaddress#1{\def\theasciiaddress{#1}}

\def\asciikeywords#1{\def\theasciikeywords{#1}}

\let\\\par\let\thelognumber\relax\let\thevolumenumber\relax
\let\thepapernumber\relax\let\thevolumeyear\relax\let\startpage\relax
\let\finishpage\relax\let\publishdate\relax\let\receiveddate\relax
\let\reviseddate\relax\let\accepteddate\relax\let\theasciititle\relax
\let\theasciiauthors\relax\let\theasciiaddress\relax
\let\theasciiabstract\relax\let\theasciikeywords\relax

\let\theasciiemail\relax

\ifplaintex
\font\logobig=cmssbx10 scaled 3836
\font\logomed=cmssbx10 scaled 2557
\else
\font\logobig=cmssbx10 scaled 4200
\font\logomed=cmssbx10 scaled 2800
\fi

\long\def\makeagttitle{   %%% start of definition of \makeagttitle
\count0=\startpage
\agt\hfill      %   Journal title (top left) 
\hbox to 45truept{\vbox to 0pt{\vglue -13truept{\logomed A\kern -.37em{\logobig 
T}\kern -.38em G}\vss}\hss}
\break
{\small Volume \thevolumenumber\ (\thevolumeyear)
\startpage--\finishpage\nl
Published: \publishdate}

\vglue .25truein

{\parskip=0pt\leftskip 0pt plus
1fil\def\\{\par\smallskip}{\Large\bf\thetitle}\par\medskip} \vglue
0.05truein

{\parskip=0pt\leftskip 0pt plus 1fil\def\\{\par}{\sc\theauthors}
\par\medskip}%
 
\vglue 0.03truein 

{\small\leftskip 25truept\rightskip 25truept{\bf Abstract}\stdspace\theabstract

{\bf AMS Classification}\stdspace\theprimaryclass
\ifx\thesecondaryclass\relax\else; \thesecondaryclass\fi\par
{\bf Keywords}\stdspace \thekeywords\par}\vglue 7truept

}   %%%% end of definition of \makeagttitle

\ifplaintex
\font\phead=cmsl9 scaled 950
\font\pnum=cmbx10 scaled 913
\font\pfoot=cmsl9 scaled 950
\headline{\vbox to 0pt{\vskip -4.5mm\line{\small\phead\ifnum
\count0=\startpage ISSN 1472-2739 (on-line) 1472-2747 (printed)
\hfill {\pnum\folio}\else\ifodd\count0\def\\{ }% 
\ifx\theshorttitle\relax\thetitle\else\theshorttitle\fi\hfill{\pnum\folio}
\else\def\\{ and }{\pnum\folio}\hfill\ifx\theshortauthors\relax\theauthors
\else\theshortauthors\fi\fi\fi}\vss}}
\footline{\vbox to 0pt{\vglue 0mm\line{\small\pfoot\ifnum\count0=\startpage
\copyright\ \gtp\hfill\else
\agt, Volume \thevolumenumber\ (\thevolumeyear)\hfill\fi}\vss}}
\else
\font=cmsl9 scaled 1050
\font\lnum=cmbx10 
\font=cmsl9 scaled 1050
\makeatletter
\def\@oddhead{{\small\lhead\ifnum\count0=\startpage ISSN 1472-2739 
(on-line) 1472-2747 (printed)\hfill {\lnum\number\count0}\else\ifodd\count0
\def\\{ }\ifx\theshorttitle\relax \thetitle \else\theshorttitle\fi\hfill
{\lnum\number\count0}\else\def\\{ and }{\lnum\number\count0}
\hfill\ifx\theshortauthors\relax 
\theauthors\else\theshortauthors\fi\fi\fi}}\def\@evenhead{\@oddhead}
\def\@oddfoot{\small\lfoot\ifnum\count0=\startpage\copyright\ \gtp\hfill\else
\agt, Volume \thevolumenumber\ (\thevolumeyear)\hfill\fi}
\def\@evenfoot{\@oddfoot}
\makeatother
\fi
\let\maketitlepage\makeagttitle
\let\makeshorttitle\maketitlepage
\let\maketitle\maketitlepage

\newwrite\gtoutfile
\long\gdef\makeheadfile{  %%% start of definition of \makeheadfile
{\def\\{, }\def\s{ }
\immediate\openout\gtoutfile head.xxx
\immediate\write\gtoutfile{To: math@arxiv.org}
\immediate\write\gtoutfile{Subject: put OR rep NNNNN:ppppp}
\immediate\write\gtoutfile{--text follows this line--}
\immediate\write\gtoutfile{Proxy-for: \ifx\theasciiauthors\relax
\theauthors\else\theasciiauthors\fi\s<\ifx\theasciiemail\relax\theemail\else\theasciiemail\fi>}
\immediate\write\gtoutfile{\noexpand\\}
\immediate\write\gtoutfile{Authors: \ifx\theasciiauthors\relax
\theauthors\else\theasciiauthors\fi}
{\def\\{ }\immediate\write\gtoutfile{Title: \ifx\theasciititle\relax
\thetitle\else\theasciititle\fi}}
\immediate\write\gtoutfile{Subj-class: GT or SG, GR etc}
\immediate\write\gtoutfile{MSC-class: \theprimaryclass\ifx\thesecondaryclass\relax\else, \thesecondaryclass\fi}
\immediate\write\gtoutfile{Journal-ref: Algebr. Geom. Topol. \thevolumenumber\s
(\thevolumeyear) \startpage-\finishpage}
\immediate\write\gtoutfile{Comments: Published by Algebraic and
Geometric Topology at}
\immediate\write\gtoutfile{\s\s\s  http://www.maths.warwick.ac.uk/agt/AGTVol\thevolumenumber/agt-\thevolumenumber-\thepapernumber.abs.html}
\immediate\write\gtoutfile{\noexpand\\}
\immediate\write\gtoutfile{}
\ifx\theasciiabstract\relax
\immediate\write\gtoutfile{\theabstract}\else
\immediate\write\gtoutfile{\theasciiabstract}\fi
\immediate\write\gtoutfile{}
\immediate\write\gtoutfile{\noexpand\\}
\immediate\write\gtoutfile{}
\immediate\closeout\gtoutfile}}  %%% end of definition of \makeheadfile

\def\maketitlepage{\makeagttitle\makeheadfile}
\let\makeshorttitle\maketitlepage
\let\maketitle\maketitlepage

\def\ifplaintex{\expandafter\ifx\csname documentclass\endcsname\relax}

\def\gtp{{\mathsurround=0pt\it $\cal G\mskip-2mu$eometry \&\ 
$\cal T\!\!$opology $\cal P\!$ublications}}  % GT publications

\def\recd{{\small Received:\qua\receiveddate\ifx\reviseddate\relax
\else\qquad Revised:\qua\reviseddate\fi\par}} 

\def\lognumber#1{\def\thelognumber{#1}}
\def\volumenumber#1{\def\thevolumenumber{#1}}
\def\volumeyear#1{\def\thevolumeyear{#1}}
\def\papernumber#1{\def\thepapernumber{#1}}
\def\pagenumbers#1#2{\def\startpage{#1}\def\finishpage{#2}}
\def\published#1{\def\publishdate{#1}}

\def\received#1{\def\receiveddate{#1}}
\def\revised#1{\def\reviseddate{#1}}
\def\accepted#1{\def\accepteddate{#1}}

\def\asciiaddress#1{\def\theasciiaddress{#1}}

\def\asciikeywords#1{\def\theasciikeywords{#1}}

\let\\\par\let\thelognumber\relax\let\thevolumenumber\relax
\let\thepapernumber\relax\let\thevolumeyear\relax\let\startpage\relax
\let\finishpage\relax\let\publishdate\relax\let\receiveddate\relax
\let\reviseddate\relax\let\accepteddate\relax\let\theasciititle\relax
\let\theasciiauthors\relax\let\theasciiaddress\relax
\let\theasciiabstract\relax\let\theasciikeywords\relax

\let\theasciiemail\relax

\ifplaintex
\font\logobig=cmssbx10 scaled 3836
\font\logomed=cmssbx10 scaled 2557
\else
\font\logobig=cmssbx10 scaled 4200
\font\logomed=cmssbx10 scaled 2800
\fi

\long\def\makeagttitle{   %%% start of definition of \makeagttitle
\count0=\startpage
\agt\hfill      %   Journal title (top left) 
\hbox to 45truept{\vbox to 0pt{\vglue -13truept{\logomed A\kern -.37em{\logobig 
T}\kern -.38em G}\vss}\hss}
\break
{\small Volume \thevolumenumber\ (\thevolumeyear)
\startpage--\finishpage\nl
Published: \publishdate}

\vglue .25truein

{\parskip=0pt\leftskip 0pt plus
1fil\def\\{\par\smallskip}{\Large\bf\thetitle}\par\medskip} \vglue
0.05truein

{\parskip=0pt\leftskip 0pt plus 1fil\def\\{\par}{\sc\theauthors}
\par\medskip}%
 
\vglue 0.03truein 

{\small\leftskip 25truept\rightskip 25truept{\bf Abstract}\stdspace\theabstract

{\bf AMS Classification}\stdspace\theprimaryclass
\ifx\thesecondaryclass\relax\else; \thesecondaryclass\fi\par
{\bf Keywords}\stdspace \thekeywords\par}\vglue 7truept

}   %%%% end of definition of \makeagttitle

\ifplaintex
\font\phead=cmsl9 scaled 950
\font\pnum=cmbx10 scaled 913
\font\pfoot=cmsl9 scaled 950
\headline{\vbox to 0pt{\vskip -4.5mm\line{\small\phead\ifnum
\count0=\startpage ISSN 1472-2739 (on-line) 1472-2747 (printed)
\hfill {\pnum\folio}\else\ifodd\count0\def\\{ }% 
\ifx\theshorttitle\relax\thetitle\else\theshorttitle\fi\hfill{\pnum\folio}
\else\def\\{ and }{\pnum\folio}\hfill\ifx\theshortauthors\relax\theauthors
\else\theshortauthors\fi\fi\fi}\vss}}
\footline{\vbox to 0pt{\vglue 0mm\line{\small\pfoot\ifnum\count0=\startpage
\copyright\ \gtp\hfill\else
\agt, Volume \thevolumenumber\ (\thevolumeyear)\hfill\fi}\vss}}
\else
\font=cmsl9 scaled 1050
\font\lnum=cmbx10 
\font=cmsl9 scaled 1050
\makeatletter
\def\@oddhead{{\small\lhead\ifnum\count0=\startpage ISSN 1472-2739 
(on-line) 1472-2747 (printed)\hfill {\lnum\number\count0}\else\ifodd\count0
\def\\{ }\ifx\theshorttitle\relax \thetitle \else\theshorttitle\fi\hfill
{\lnum\number\count0}\else\def\\{ and }{\lnum\number\count0}
\hfill\ifx\theshortauthors\relax 
\theauthors\else\theshortauthors\fi\fi\fi}}\def\@evenhead{\@oddhead}
\def\@oddfoot{\small\lfoot\ifnum\count0=\startpage\copyright\ \gtp\hfill\else
\agt, Volume \thevolumenumber\ (\thevolumeyear)\hfill\fi}
\def\@evenfoot{\@oddfoot}
\makeatother
\fi
\let\maketitlepage\makeagttitle
\let\makeshorttitle\maketitlepage
\let\maketitle\maketitlepage

\newwrite\gtoutfile
\long\gdef\makeheadfile{  %%% start of definition of \makeheadfile
{\def\\{, }\def\s{ }
\immediate\openout\gtoutfile head.xxx
\immediate\write\gtoutfile{To: math@arxiv.org}
\immediate\write\gtoutfile{Subject: put OR rep NNNNN:ppppp}
\immediate\write\gtoutfile{--text follows this line--}
\immediate\write\gtoutfile{Proxy-for: \ifx\theasciiauthors\relax
\theauthors\else\theasciiauthors\fi\s<\ifx\theasciiemail\relax\theemail\else\theasciiemail\fi>}
\immediate\write\gtoutfile{\noexpand\\}
\immediate\write\gtoutfile{Authors: \ifx\theasciiauthors\relax
\theauthors\else\theasciiauthors\fi}
{\def\\{ }\immediate\write\gtoutfile{Title: \ifx\theasciititle\relax
\thetitle\else\theasciititle\fi}}
\immediate\write\gtoutfile{Subj-class: GT or SG, GR etc}
\immediate\write\gtoutfile{MSC-class: \theprimaryclass\ifx\thesecondaryclass\relax\else, \thesecondaryclass\fi}
\immediate\write\gtoutfile{Journal-ref: Algebr. Geom. Topol. \thevolumenumber\s
(\thevolumeyear) \startpage-\finishpage}
\immediate\write\gtoutfile{Comments: Published by Algebraic and
Geometric Topology at}
\immediate\write\gtoutfile{\s\s\s  http://www.maths.warwick.ac.uk/agt/AGTVol\thevolumenumber/agt-\thevolumenumber-\thepapernumber.abs.html}
\immediate\write\gtoutfile{\noexpand\\}
\immediate\write\gtoutfile{}
\ifx\theasciiabstract\relax
\immediate\write\gtoutfile{\theabstract}\else
\immediate\write\gtoutfile{\theasciiabstract}\fi
\immediate\write\gtoutfile{}
\immediate\write\gtoutfile{\noexpand\\}
\immediate\write\gtoutfile{}
\immediate\closeout\gtoutfile}}  %%% end of definition of \makeheadfile

\def\maketitlepage{\makeagttitle\makeheadfile}
\let\makeshorttitle\maketitlepage
\let\maketitle\maketitlepage

\def\ifplaintex{\expandafter\ifx\csname documentclass\endcsname\relax}

\def\gtp{{\mathsurround=0pt\it $\cal G\mskip-2mu$eometry \&\ 
$\cal T\!\!$opology $\cal P\!$ublications}}  % GT publications

\def\recd{{\small Received:\qua\receiveddate\ifx\reviseddate\relax
\else\qquad Revised:\qua\reviseddate\fi\par}} 

\def\lognumber#1{\def\thelognumber{#1}}
\def\volumenumber#1{\def\thevolumenumber{#1}}
\def\volumeyear#1{\def\thevolumeyear{#1}}
\def\papernumber#1{\def\thepapernumber{#1}}
\def\pagenumbers#1#2{\def\startpage{#1}\def\finishpage{#2}}
\def\published#1{\def\publishdate{#1}}

\def\received#1{\def\receiveddate{#1}}
\def\revised#1{\def\reviseddate{#1}}
\def\accepted#1{\def\accepteddate{#1}}

\def\asciiaddress#1{\def\theasciiaddress{#1}}

\def\asciikeywords#1{\def\theasciikeywords{#1}}

\let\\\par\let\thelognumber\relax\let\thevolumenumber\relax
\let\thepapernumber\relax\let\thevolumeyear\relax\let\startpage\relax
\let\finishpage\relax\let\publishdate\relax\let\receiveddate\relax
\let\reviseddate\relax\let\accepteddate\relax\let\theasciititle\relax
\let\theasciiauthors\relax\let\theasciiaddress\relax
\let\theasciiabstract\relax\let\theasciikeywords\relax

\let\theasciiemail\relax

\ifplaintex
\font\logobig=cmssbx10 scaled 3836
\font\logomed=cmssbx10 scaled 2557
\else
\font\logobig=cmssbx10 scaled 4200
\font\logomed=cmssbx10 scaled 2800
\fi

\long\def\makeagttitle{   %%% start of definition of \makeagttitle
\count0=\startpage
\agt\hfill      %   Journal title (top left) 
\hbox to 45truept{\vbox to 0pt{\vglue -13truept{\logomed A\kern -.37em{\logobig 
T}\kern -.38em G}\vss}\hss}
\break
{\small Volume \thevolumenumber\ (\thevolumeyear)
\startpage--\finishpage\nl
Published: \publishdate}

\vglue .25truein

{\parskip=0pt\leftskip 0pt plus
1fil\def\\{\par\smallskip}{\Large\bf\thetitle}\par\medskip} \vglue
0.05truein

{\parskip=0pt\leftskip 0pt plus 1fil\def\\{\par}{\sc\theauthors}
\par\medskip}%
 
\vglue 0.03truein 

{\small\leftskip 25truept\rightskip 25truept{\bf Abstract}\stdspace\theabstract

{\bf AMS Classification}\stdspace\theprimaryclass
\ifx\thesecondaryclass\relax\else; \thesecondaryclass\fi\par
{\bf Keywords}\stdspace \thekeywords\par}\vglue 7truept

}   %%%% end of definition of \makeagttitle

\ifplaintex
\font\phead=cmsl9 scaled 950
\font\pnum=cmbx10 scaled 913
\font\pfoot=cmsl9 scaled 950
\headline{\vbox to 0pt{\vskip -4.5mm\line{\small\phead\ifnum
\count0=\startpage ISSN 1472-2739 (on-line) 1472-2747 (printed)
\hfill {\pnum\folio}\else\ifodd\count0\def\\{ }% 
\ifx\theshorttitle\relax\thetitle\else\theshorttitle\fi\hfill{\pnum\folio}
\else\def\\{ and }{\pnum\folio}\hfill\ifx\theshortauthors\relax\theauthors
\else\theshortauthors\fi\fi\fi}\vss}}
\footline{\vbox to 0pt{\vglue 0mm\line{\small\pfoot\ifnum\count0=\startpage
\copyright\ \gtp\hfill\else
\agt, Volume \thevolumenumber\ (\thevolumeyear)\hfill\fi}\vss}}
\else
\font=cmsl9 scaled 1050
\font\lnum=cmbx10 
\font=cmsl9 scaled 1050
\makeatletter
\def\@oddhead{{\small\lhead\ifnum\count0=\startpage ISSN 1472-2739 
(on-line) 1472-2747 (printed)\hfill {\lnum\number\count0}\else\ifodd\count0
\def\\{ }\ifx\theshorttitle\relax \thetitle \else\theshorttitle\fi\hfill
{\lnum\number\count0}\else\def\\{ and }{\lnum\number\count0}
\hfill\ifx\theshortauthors\relax 
\theauthors\else\theshortauthors\fi\fi\fi}}\def\@evenhead{\@oddhead}
\def\@oddfoot{\small\lfoot\ifnum\count0=\startpage\copyright\ \gtp\hfill\else
\agt, Volume \thevolumenumber\ (\thevolumeyear)\hfill\fi}
\def\@evenfoot{\@oddfoot}
\makeatother
\fi
\let\maketitlepage\makeagttitle
\let\makeshorttitle\maketitlepage
\let\maketitle\maketitlepage

\newwrite\gtoutfile
\long\gdef\makeheadfile{  %%% start of definition of \makeheadfile
{\def\\{, }\def\s{ }
\immediate\openout\gtoutfile head.xxx
\immediate\write\gtoutfile{To: math@arxiv.org}
\immediate\write\gtoutfile{Subject: put OR rep NNNNN:ppppp}
\immediate\write\gtoutfile{--text follows this line--}
\immediate\write\gtoutfile{Proxy-for: \ifx\theasciiauthors\relax
\theauthors\else\theasciiauthors\fi\s<\ifx\theasciiemail\relax\theemail\else\theasciiemail\fi>}
\immediate\write\gtoutfile{\noexpand\\}
\immediate\write\gtoutfile{Authors: \ifx\theasciiauthors\relax
\theauthors\else\theasciiauthors\fi}
{\def\\{ }\immediate\write\gtoutfile{Title: \ifx\theasciititle\relax
\thetitle\else\theasciititle\fi}}
\immediate\write\gtoutfile{Subj-class: GT or SG, GR etc}
\immediate\write\gtoutfile{MSC-class: \theprimaryclass\ifx\thesecondaryclass\relax\else, \thesecondaryclass\fi}
\immediate\write\gtoutfile{Journal-ref: Algebr. Geom. Topol. \thevolumenumber\s
(\thevolumeyear) \startpage-\finishpage}
\immediate\write\gtoutfile{Comments: Published by Algebraic and
Geometric Topology at}
\immediate\write\gtoutfile{\s\s\s  http://www.maths.warwick.ac.uk/agt/AGTVol\thevolumenumber/agt-\thevolumenumber-\thepapernumber.abs.html}
\immediate\write\gtoutfile{\noexpand\\}
\immediate\write\gtoutfile{}
\ifx\theasciiabstract\relax
\immediate\write\gtoutfile{\theabstract}\else
\immediate\write\gtoutfile{\theasciiabstract}\fi
\immediate\write\gtoutfile{}
\immediate\write\gtoutfile{\noexpand\\}
\immediate\write\gtoutfile{}
\immediate\closeout\gtoutfile}}  %%% end of definition of \makeheadfile

\def\maketitlepage{\makeagttitle\makeheadfile}
\let\makeshorttitle\maketitlepage
\let\maketitle\maketitlepage

\def\ifplaintex{\expandafter\ifx\csname documentclass\endcsname\relax}

\def\gtp{{\mathsurround=0pt\it $\cal G\mskip-2mu$eometry \&\ 
$\cal T\!\!$opology $\cal P\!$ublications}}  % GT publications

\def\recd{{\small Received:\qua\receiveddate\ifx\reviseddate\relax
\else\qquad Revised:\qua\reviseddate\fi\par}} 

\def\lognumber#1{\def\thelognumber{#1}}
\def\volumenumber#1{\def\thevolumenumber{#1}}
\def\volumeyear#1{\def\thevolumeyear{#1}}
\def\papernumber#1{\def\thepapernumber{#1}}
\def\pagenumbers#1#2{\def\startpage{#1}\def\finishpage{#2}}
\def\published#1{\def\publishdate{#1}}

\def\received#1{\def\receiveddate{#1}}
\def\revised#1{\def\reviseddate{#1}}
\def\accepted#1{\def\accepteddate{#1}}

\def\asciiaddress#1{\def\theasciiaddress{#1}}

\def\asciikeywords#1{\def\theasciikeywords{#1}}

\let\\\par\let\thelognumber\relax\let\thevolumenumber\relax
\let\thepapernumber\relax\let\thevolumeyear\relax\let\startpage\relax
\let\finishpage\relax\let\publishdate\relax\let\receiveddate\relax
\let\reviseddate\relax\let\accepteddate\relax\let\theasciititle\relax
\let\theasciiauthors\relax\let\theasciiaddress\relax
\let\theasciiabstract\relax\let\theasciikeywords\relax

\let\theasciiemail\relax

\ifplaintex
\font\logobig=cmssbx10 scaled 3836
\font\logomed=cmssbx10 scaled 2557
\else
\font\logobig=cmssbx10 scaled 4200
\font\logomed=cmssbx10 scaled 2800
\fi

\long\def\makeagttitle{   %%% start of definition of \makeagttitle
\count0=\startpage
\agt\hfill      %   Journal title (top left) 
\hbox to 45truept{\vbox to 0pt{\vglue -13truept{\logomed A\kern -.37em{\logobig 
T}\kern -.38em G}\vss}\hss}
\break
{\small Volume \thevolumenumber\ (\thevolumeyear)
\startpage--\finishpage\nl
Published: \publishdate}

\vglue .25truein

{\parskip=0pt\leftskip 0pt plus
1fil\def\\{\par\smallskip}{\Large\bf\thetitle}\par\medskip} \vglue
0.05truein

{\parskip=0pt\leftskip 0pt plus 1fil\def\\{\par}{\sc\theauthors}
\par\medskip}%
 
\vglue 0.03truein 

{\small\leftskip 25truept\rightskip 25truept{\bf Abstract}\stdspace\theabstract

{\bf AMS Classification}\stdspace\theprimaryclass
\ifx\thesecondaryclass\relax\else; \thesecondaryclass\fi\par
{\bf Keywords}\stdspace \thekeywords\par}\vglue 7truept

}   %%%% end of definition of \makeagttitle

\ifplaintex
\font\phead=cmsl9 scaled 950
\font\pnum=cmbx10 scaled 913
\font\pfoot=cmsl9 scaled 950
\headline{\vbox to 0pt{\vskip -4.5mm\line{\small\phead\ifnum
\count0=\startpage ISSN 1472-2739 (on-line) 1472-2747 (printed)
\hfill {\pnum\folio}\else\ifodd\count0\def\\{ }% 
\ifx\theshorttitle\relax\thetitle\else\theshorttitle\fi\hfill{\pnum\folio}
\else\def\\{ and }{\pnum\folio}\hfill\ifx\theshortauthors\relax\theauthors
\else\theshortauthors\fi\fi\fi}\vss}}
\footline{\vbox to 0pt{\vglue 0mm\line{\small\pfoot\ifnum\count0=\startpage
\copyright\ \gtp\hfill\else
\agt, Volume \thevolumenumber\ (\thevolumeyear)\hfill\fi}\vss}}
\else
\font=cmsl9 scaled 1050
\font\lnum=cmbx10 
\font=cmsl9 scaled 1050
\makeatletter
\def\@oddhead{{\small\lhead\ifnum\count0=\startpage ISSN 1472-2739 
(on-line) 1472-2747 (printed)\hfill {\lnum\number\count0}\else\ifodd\count0
\def\\{ }\ifx\theshorttitle\relax \thetitle \else\theshorttitle\fi\hfill
{\lnum\number\count0}\else\def\\{ and }{\lnum\number\count0}
\hfill\ifx\theshortauthors\relax 
\theauthors\else\theshortauthors\fi\fi\fi}}\def\@evenhead{\@oddhead}
\def\@oddfoot{\small\lfoot\ifnum\count0=\startpage\copyright\ \gtp\hfill\else
\agt, Volume \thevolumenumber\ (\thevolumeyear)\hfill\fi}
\def\@evenfoot{\@oddfoot}
\makeatother
\fi
\let\maketitlepage\makeagttitle
\let\makeshorttitle\maketitlepage
\let\maketitle\maketitlepage

\newwrite\gtoutfile
\long\gdef\makeheadfile{  %%% start of definition of \makeheadfile
{\def\\{, }\def\s{ }
\immediate\openout\gtoutfile head.xxx
\immediate\write\gtoutfile{To: math@arxiv.org}
\immediate\write\gtoutfile{Subject: put OR rep NNNNN:ppppp}
\immediate\write\gtoutfile{--text follows this line--}
\immediate\write\gtoutfile{Proxy-for: \ifx\theasciiauthors\relax
\theauthors\else\theasciiauthors\fi\s<\ifx\theasciiemail\relax\theemail\else\theasciiemail\fi>}
\immediate\write\gtoutfile{\noexpand\\}
\immediate\write\gtoutfile{Authors: \ifx\theasciiauthors\relax
\theauthors\else\theasciiauthors\fi}
{\def\\{ }\immediate\write\gtoutfile{Title: \ifx\theasciititle\relax
\thetitle\else\theasciititle\fi}}
\immediate\write\gtoutfile{Subj-class: GT or SG, GR etc}
\immediate\write\gtoutfile{MSC-class: \theprimaryclass\ifx\thesecondaryclass\relax\else, \thesecondaryclass\fi}
\immediate\write\gtoutfile{Journal-ref: Algebr. Geom. Topol. \thevolumenumber\s
(\thevolumeyear) \startpage-\finishpage}
\immediate\write\gtoutfile{Comments: Published by Algebraic and
Geometric Topology at}
\immediate\write\gtoutfile{\s\s\s  http://www.maths.warwick.ac.uk/agt/AGTVol\thevolumenumber/agt-\thevolumenumber-\thepapernumber.abs.html}
\immediate\write\gtoutfile{\noexpand\\}
\immediate\write\gtoutfile{}
\ifx\theasciiabstract\relax
\immediate\write\gtoutfile{\theabstract}\else
\immediate\write\gtoutfile{\theasciiabstract}\fi
\immediate\write\gtoutfile{}
\immediate\write\gtoutfile{\noexpand\\}
\immediate\write\gtoutfile{}
\immediate\closeout\gtoutfile}}  %%% end of definition of \makeheadfile

\def\maketitlepage{\makeagttitle\makeheadfile}
\let\makeshorttitle\maketitlepage
\let\maketitle\maketitlepage

\lognumber{1}
\volumenumber{2}
\volumeyear{2002}
\papernumber{1}
\published{12 January 2002}
\pagenumbers{1}{9}
\received{22 August 2001}
\revised{8 January 2002}
\accepted{10 January 2002}

\usepackage{amssymb}
\usepackage{amsmath,amscd}

\newcommand{\N}{{\mathbb N}}

\newcommand{\Mon}{\mathop{\mathrm{Mon}}}

\newtheorem{theorem}{Theorem}
\newtheorem{lemma}[theorem]{Lemma}

\newtheorem{corollary}[theorem]{Corollary}

\theoremstyle{definition}
\newtheorem{definition}[theorem]{Definition}

\title{Bihomogeneity of solenoids}

\author{Alex Clark\\Robbert Fokkink}

\address{University of North Texas, Department of Mathematics\\Denton
TX 76203-1430, U.S.A.}
\secondaddress{Technische Universiteit Delft,
Faculty of Information Technology and Systems\\Division Mediamatica,
P.O. Box 5031, 2600 GA Delft, Netherlands}

\asciiaddress{University of North Texas, Department of Mathematics\\Denton
TX 76203-1430, USA\\and\\Technische Universiteit Delft,
Faculty of Information Technology and Systems\\Division Mediamatica,
P.O. Box 5031, 2600 GA Delft, Netherlands}

\email{alexc@unt.edu, r.j.fokkink@its.tudelft.nl}

\begin{document}

\begin{abstract}
Solenoids are inverse limit spaces over regular covering maps of closed
manifolds.
M.C. McCord has shown that solenoids are topologically homogeneous and
that they are principal bundles with a profinite structure group.
We show that if a solenoid is bihomogeneous,
then its structure group
contains an open abelian subgroup. This leads
to new
examples of homogeneous continua that are not bihomogeneous.
\end{abstract}

\primaryclass{54F15}
\secondaryclass{55R10}
\keywords{Homogeneous continuum, covering space, profinite group,\break
principal bundle}
\asciikeywords{Homogeneous continuum, covering space, profinite group,
principal bundle}

\makeshorttitle

A topological space $X$ is \emph{homogeneous\/} if for every pair of points
$x,y\in X$ there is a homeomorphism $h:X\rightarrow X$ satisfying
$h\left(
x\right) =y$. The space is \emph{bihomogeneous\/} if for each
such pair there is a homeomorphism satisfying $h\left( x\right) =y$ and
$
h\left( y\right) =x.$
A compact and connected space is called a \emph{continuum}.
Knaster and Van Dantzig asked whether a
homogeneous
continuum is necessarily bihomogeneous. This was settled in the negative
by
Krystyna Kuperberg \cite{Kuperberg}. Subsequent counterexamples were
given
by Minc, Kawamura and Greg Kuperberg \cite{Minc, Kawamura, GKuperberg}.
The
counterexamples in \cite{Kuperberg,GKuperberg} are locally connected.
Ungar \cite{Ungar} has studied stronger types of homogeneity
conditions and showed that these conditions imply local connectivity.

A solenoid $M_{\infty }$ is an inverse limit space
over closed connected manifolds with bonding maps that are covering maps.
We shall silently
assume that the bonding maps are not $1-1$, so that
$M_{\infty }$ is not locally connected.
McCord~\cite{McCord} has shown that solenoids are homogeneous
provided that compositions of the bonding covering maps are regular. Minc~\cite{Minc}
presented an example of a homogeneous but not bihomogeneous infinite-dimensional 
continuum similar to a solenoid, and Krystyna Kuperberg~\cite{K2} 
observed that a similar construction could be used to construct a 
finite-dimensional solenoid which is homogeneous but not bihomogeneous.
We shall show
that $M_{\infty }$ is bihomogeneous only if a certain
condition
related to commutativity (or lack thereof) of $\pi _{1}\left(
M_{i}\right) $
is met. In case the solenoid
is $2$-dimensional, the condition is both necessary and sufficient.

\section{Path-components of solenoids as left-cosets of the structure
group}

A {\it (strong) solenoid} $M_\infty$ is an inverse-limit space of closed
manifolds
$M_i$ with bonding maps $p_i\colon M_{i+1}\to M_i$ for $i\in\N$ which
are covering maps, such that any composition $p_{i+k}\circ \ldots\circ
p_i$
is regular.
Solenoids are homogeneous spaces
and they have dense path-components.
\smallbreak
A $G$-bundle $(E,B,p,F)$ is {\it principal} if the
structure group $G$ acts effectively on the fibers.
As a consequence, the fiber
$F$ is homeomorphic to $G$, and
$G$ is naturally equivalent to
the group of deck-transformations.
\begin{theorem}[McCord, \protect\cite{McCord}]
Suppose that $M_\infty=\lim_{\leftarrow}(M_i,p_i)$ is a solenoid.
Let $\pi_0\colon M_\infty\to M_0$ be the projection onto the first
coordinate
and let $\Gamma_0=\pi_0^{-1}(m_0)$ be a fiber.
Then $(M_{\infty },M_{0},\pi_0,\Gamma _{0})$ is a principal-bundle.
\end{theorem}
The projection $\pi_0$ is not to be confused with a homotopy group.
Note that a solenoid $\lim_{\leftarrow }(M_{i},f_{i})$ is a principal
bundle over any $M_i$ and we have singled out $M_0$.
The spaces $M_i$ are called the {\it factor spaces} of the solenoid.
We think of the fundamental groups $\pi_1(M_i)$ as (normal) subgroups
of $\pi_1(M_0)$.
The structure group $\Gamma_0$ is isomorphic to the profinite group
$\lim_{\leftarrow}\pi_1(M_0)/\pi_1(M_i)$.
\medbreak
Choose base-points $m_i\in M_i$ such that $p_i(m_i)=m_{i-1}$, so
$m_\infty=(m_i)$ is an element of $M_\infty$.
We identify the structure group $\Gamma_0$ with the fiber of $m_0$
and we identify $m_\infty$ with the unit element of $\Gamma_0$.
The fundamental group $\pi_1(M_0)$ acts on the base-point fiber
$\Gamma_0$ by path lifting: for $g\in\Gamma_0$ and $\gamma\in\pi_1(M_0,m_0)$,
define $g\circ\gamma$ as the end-point of the lifted path $\tilde\gamma$
starting from the initial-point $g$.
One verifies that this right action of $\pi_1(M_0)$
commutes with left multiplication of $\Gamma_0$.
More precisely,
suppose that $h$ is a deck-transformation and that $\tilde\gamma$ is a
lifted path with initial-point $g$.
Then $h(\tilde\gamma)$ has initial point $h(g)$ and end-point
$h(g\circ\gamma)$.
Identify the structure group with the group of deck-transformations,
so we get that $(hg)\circ\gamma=h(g\circ\gamma)$.

\begin{definition}
Suppose that $(M_{\infty },M_0,\pi_0,\Gamma_0)$ is a solenoid.
We shall call the $\pi_1(M_0)$-orbit of $e\in\Gamma_0$
the characteristic group and we shall denote it by $\gamma_0$.
Let $K_\infty\subset\pi_1(M_0)$ be the intersection of all
$\pi_1(M_i)$.
Then $\gamma_0$ is isomorphic to $\pi_1(M_0)/K_\infty$ and we shall
refer to $K_\infty$ as the kernel of $\pi_1(M_0)$.
\end{definition}

Our definition deviates from the common terminology, as in \cite{Steenrod},
where the equivalence class of $\gamma_0$ under inner automorphisms of
$\Gamma_0$ is called the characteristic class.
Note that $\gamma_0$ inherits a topology from $\Gamma_0$.

\begin{lemma}\label{leftcosets}
The path components of a solenoid are naturally equivalent to the left
cosets
$\Gamma_0/\gamma_0$.
\end{lemma}
\begin{proof}
Suppose that $x,y\in\Gamma_0$ are elements of the base-point fiber.
Then $x\circ\gamma=y$ for some $\gamma\in\pi_1(M_0)$ if and only
if there exists a path $\tilde\gamma\subset M_\infty$ that connects
$x$ to $y$.
\end{proof}

If we replace the base space $M_0$ by $M_i$ for some index $i$,
then we get a principal bundle $(M_\infty,M_i,\pi_i,\Gamma_i)$,
where $\Gamma_i\subset \Gamma_0$ is the subgroup of transformations that
leave $M_i$ invariant.
The topology of $\Gamma_0$ is induced by taking the $\Gamma_i$ as
an open neighborhood base of the identity.
One verifies that the charateristic group of the bundle, denoted
$\gamma_i$,
is equal to $\gamma_0\cap\Gamma_i$.
Hence the $\gamma_i$ are open subgroups of $\gamma_0$.

\begin{lemma}
For $j>i$ the inclusion $\Gamma _{j}\subset \Gamma _{i}$ induces a
natural isomorphism between $\Gamma _{j}/\gamma _{j}$
and $\Gamma _{i}/\gamma _{i}.$
\end{lemma}

Since path components are dense in $M_{\infty },$ the
characteristic subgroups $\gamma _{i}$ are dense in $\Gamma _{i}$.

\section{The permutation of path-components by self-\break hom\-eo\-morphisms}

A solenoid $M_\infty$ can be represented as a subspace of $\prod M_i$,
the Cartesian product of its factor spaces.
We identify
$M_i$ with the subspace of $\prod M_i$ defined by:
$$M_i=\{(x_j) : x_j\in M_j,\ x_j=p^i_j(x_i) \text{ if } j\leq i,\
x_j=m_j \text{ if }j>i\}$$
where $p^i_j\colon M_i\to M_j$ is a composition of bonding maps.
In this representation,
the factor spaces $M_i$ and $M_\infty$ all have the same base-point.
\medbreak
A \emph{morphism} between fiber bundles
can be represented by a commutative diagram:
\begin{equation*}
\begin{array}{rcl}
E_{1}& \overset{h}{\longrightarrow } & E_{2}\  \\
p_1 \downarrow &  &\downarrow p_2 \\
B_{1} & \overset{f}{\longrightarrow } & B_{2}\
\end{array}
\end{equation*}
We shall say that $h$ is the {\it lifted map} and that $f$ is the {\it
base-map}.
We say that morphisms are homotopic if their base-maps are.
By the unique path-lifting property, a morphism between bundles
with a totally disconnected fiber is determined by the base-map
$f\colon B_1\to B_2$
and the image under $h$ of a single element of $E_1$.
For pointed spaces, therefore, a bundle-morphism is determined by the
base-map only.
This implies that,
for principal bundles with a totally disconnected fiber,
bundle-morphisms
commute with deck-transformations; i.e., for
a lifted map $h$ and a deck-transformation $\varphi\colon E_1\to E_1$,
we have that $h\circ\varphi=\psi\circ h$ for some
deck-transformation $\psi\colon E_2\to E_2$.

\begin{lemma}\label{induceshomomorphism}
Suppose that $(E_{i},B_{i},p_{i},\Gamma_{i})$ are principal
$\Gamma_i$-bundles with a totally disconnected fiber
(for $i=1,2$).
Then a base-point preserving
bundle-morphism induces a homomorphism of the structure group.
Furthermore, homotopic morphisms induce the same homomorphism.
\end{lemma}
\begin{proof}
First note that the lifted map $h$ maps $\Gamma_1$ to $\Gamma_2$.
Deck-transformations are (left) translations
$x\to ax$ of the
base-point fiber $\Gamma_i$ ($i=1,2$).
Since a bundle-morphism commutes with deck-transformations,
$h\colon \Gamma_1\to \Gamma_2$ satisfies $h(ax)=f(a)h(x)$ for some
$f\colon \Gamma_1\to \Gamma_2$.
Substitute $x=e$ to find that $h(ax)=h(a)h(x)$.
Now homotopic bundle-morphisms give homotopic homomorphisms $h\colon
\Gamma_1\to\Gamma_2$.
Since the groups are totally disconnected, the homomorphisms are
necessarily
the same.
\end{proof}

We shall say that a bundle morphism of a solenoid is an
\textit{automorphism}
if the commutative diagram can be extended on the right-hand side
\begin{equation*}
\begin{array}{rclrcl}
M_{\infty}& \overset{h_1}{\longrightarrow } & M_{\infty}&
\overset{h_2}{\longrightarrow } & M_{\infty}\\
\pi_j \downarrow &  &\downarrow \pi_i &  &\downarrow \pi_k\\
M_{j} & \overset{f_1}{\longrightarrow } & M_{i}&
\overset{f_2}{\longrightarrow } & M_{k}
\end{array}
\end{equation*}

such that $f_2\circ f_1$ is homotopic to $p^j_k$. We shall say that
$h_1$ is the
inverse of $h_2$.
For instance, the covering
projection $p^j_i\colon M_j\to M_i$ with
lifted map $id_{M_\infty}$
yields an automorphism.
We show that for every self-homeomorphism of a
solenoid,
there is an automorphism that acts in the same way on the space of
path-components.

\begin{theorem}\label{homotopicmorphism}
Suppose that $h$ is a base-point preserving
self-homeomorphism of a solenoid $M_{\infty }$.
Then $h$ is homotopic to the lifted map of
an automorphism of $M_\infty$.
\end{theorem}

\begin{proof}
Since $M_0$ is an ANR, the composition
$\pi_0\circ h\colon M_\infty\to M_0$
extends to $H\colon U\to M_0$ for
a neighborhood of $M_\infty\subset U$ in $\prod M_i$.
The restriction $H\colon M_i\to M_0$ is well-defined for
sufficiently large $i$.
Note that $H$ preserves the base-point of $M_i$.
For sufficiently large $i$, the maps $H\circ \pi_i$ and $\pi_0\circ h$
are homotopic.
By the homotopy lifting property, $H\circ \pi_i$ can then be lifted to
$\tilde H\colon M_\infty\to M_\infty$, which is homotopic to $h$.

Now apply the same argument to $\pi_i\circ h^{-1}$ to find
a map $G\colon M_j\to M_i$ for sufficiently large $j$ which
can be lifted to $\tilde G\colon M_\infty\to M_\infty$.
By choosing $j$ and $i$ sufficiently large, the composition
$H\circ G\colon M_j\to M_0$ gets arbitrarily close
to and hence homotopic to the covering map $p^j_0$.
\end{proof}

Theorem \ref{homotopicmorphism} and Lemma \ref{induceshomomorphism}
describe
how a self-homeomorphism acts on path-components of
a solenoid (provided that it preserves the
base-point).

\begin{lemma}
Suppose that $h$ is the lifted map of an automorphism of a solenoid
$M_\infty$.
For some index $i$, $h$ induces a monomorphism
$\hat h\colon\Gamma_i\to \Gamma_0$ such that $\hat
{h}^{-1}(\gamma_0)=\gamma_i$
and $\hat h(\Gamma_i)$ is an open subgroup of $\Gamma_0$.
\end{lemma}
\begin{proof}
By Lemma \ref{induceshomomorphism} we know that $h$ induces
a homomorphism $\hat h\colon \Gamma_i\to\Gamma_0$.
Since homeomorphisms preserve path-components,
Lemma \ref{leftcosets} implies that $h$ induces a homomorphism
$\Gamma_i/\gamma_i\to \Gamma_0/\gamma_0$.
Since $h$ is the lifted map of an automorphism, it has an inverse $g$
which induces a homomorphism $\hat{g}\colon\Gamma_j\to \Gamma_0$.
The composition $\hat g\circ\hat h$, which is defined on an open
subgroup,
is equal to
the identity.
By Lemma \ref{induceshomomorphism},
$\hat g\circ\hat h$ is equal to the homomorphism induced
by $p^j_i$, which is the identity.
\end{proof}

\section{An algebraic condition for bihomogeneity}

\begin{definition}
Suppose that $\Gamma_0$ is the structure group of a solenoid with
characteristic group $\gamma_0$.
We define $\Mon(\Gamma_0,\gamma_0)$ as the set of monomorphisms
$f\colon \Gamma_i\to \Gamma_0$, such that $f(\gamma_i)=\gamma_0\cap
f(\Gamma_i)$.
\end{definition}

We say that an element of $\Mon(\Gamma_0,\gamma_0)$ is a
{\emph{characteristic automorphism}}.
A self-homeomorphism $H$ of $M_\infty$ need not preserve
the base-point.
It can however be represented as a composition
of a homeomorphism $h$ that preserves the path-component of
the base-point and a deck-transformation.
Since $h$ is homotopic to a base-point preserving homeomorphism,
$H$ permutes the path-components in the same way as a composition of
a base-point preserving homeomorphism and a deck-transformation.
In terms of $\Gamma_0/\gamma_0$,
this is a composition of a characteristic automorphism
$\varphi$ and a left translation $z\to wz$ of $\Gamma_0$.

\begin{definition}
We say that a solenoid is \emph{algebraically bihomogeneous\/} if
it satisfies the following condition.
For every $x,y\in \Gamma _{0}$ there are elements $w\in\Gamma_0$
and $\varphi\in\Mon(\Gamma_0,\gamma_0)$ such that
$z\to w\varphi(z)$ switches the residue classes
$x\mathop{\mathrm{mod}}\gamma_0$ and
$y\mathop{\mathrm{mod}}\gamma_0$.
\end{definition}

Obviously, bihomogeneity implies algebraic bihomogeneity.
The condition of algebraic bihomogeneity may seem awkward, but
fortunately there is a simpler characterization as we shall
see below.
We denote $x\sim y$ if $x,y$ are in the same residue class of
$\gamma_{0}$.

\begin{lemma}\label{inverse}
A solenoid $M_{\infty }$ is algebraically bihomogeneous if and if only for
every $z\in \Gamma _{0}$ there is a characteristic automorphism
$\varphi$ such that $\varphi (z)\sim z^{-1}$.
\end{lemma}

\begin{proof}
Suppose that $M_{\infty }$ is algebraically bihomogeneous.
For every $z\in
\Gamma _{0}$, we can switch the cosets of $z$ and $e$.
More precisely, there exists a
$w\in \Gamma_0$ and a $\varphi\in\Mon(\Gamma_0,\gamma_0)$
such that $zg=w\varphi (e)$ and $eg^{\prime }=w\varphi (z)$
for $g,g^{\prime }\in \gamma _{0}$.
Since $\varphi (e)=e,$ it follows that
$w=zg$ and $\varphi (z)=g^{-1}z^{-1}g^{\prime }$. Compose $
\varphi $ with the inner automorphism $x\rightarrow gxg^{-1}$ to obtain
$
\psi \in \Mon(\Gamma _{0},\gamma_0)$ satisfying $\psi
\left( z\right) \sim z^{-1}.$

If $\varphi (z)=z^{-1}g$ for some $g\in \gamma _{0}$, then
compose $\varphi $ with the inner automorphism $x\rightarrow gxg^{-1}$
to
get $\psi \in \Mon(\Gamma _{0},\gamma_0)$ satisfying
$\psi (z)=gz^{-1}$. Then $x\rightarrow zg^{-1}\psi (x)$ switches the
cosets of~$e$ and~$z$.
This implies algebraic bihomogeneity.
\end{proof}

Since $
z\rightarrow z^{-1}$ is a homomorphism if and only if the group is abelian, we have the
following corollary.

\begin{corollary}
A solenoid with an abelian structure group $\Gamma _{i}$ is
algebraically
bihomogeneous.
\end{corollary}

This condition is automatically met if $\pi _{1}\left(
M_{i}\right) $ is abelian.

\begin{lemma}
Suppose that $\gamma _{0}$ is a characteristic group. Then
$\Mon(\Gamma _{0},\gamma_0)$ is countable.
\end{lemma}

\begin{proof}
There are only countably many subgroups $\gamma_i$ and each of these
is finitely generated.
Hence, there are only finitely many homomorphisms
$f\colon\gamma_i\to\gamma_0$.
Since characteristic automorphisms are determined by their action on
some $\gamma_i$, the result follows.
\end{proof}

\begin{theorem}
\label{main}Let $M_{\infty }$ be a bihomogeneous solenoid with structure
group $\Gamma _{0}$. Then $\Gamma _{0}$ contains an open abelian
subgroup.
\end{theorem}

\begin{proof}
Suppose that $\varphi \colon \Gamma _{j}\rightarrow \Gamma _{0}$ is a
characteristic automorphism. For $g\in \gamma _{0}$ define the subset $
V(\varphi ,g)=\{z\in \Gamma _{j}\colon z\varphi (z)=g\}\subset \Gamma
_{0}.$
As $\varphi $ ranges over $\Mon(\Gamma_{0},\gamma_0)$
and $g$ ranges over $\gamma _{0}$, the countable family of all
$V(\varphi ,g)
$ covers $\Gamma _{0}$ by Lemma \ref{inverse}.
Hence one of
these sets, say $V(\varphi _{0},g_{0})$, is of second category in
$\Gamma
_{0}$.
It follows that $K=\{z\in \Gamma _{0}\colon
z\varphi
_{0}(z)=g_{0}\}$ is closed with non-empty interior in $\Gamma _{0}.$
Since $K$ has non-empty interior,
there exist a $z_{0}\in K$ and a neighborhood $V$ of $e$ such that
$\varphi
(z_{0}\delta )=\delta ^{-1}z_{0}^{-1}g_{0}$ for all $\delta \in V$. It
follows that $\varphi (\delta )=g_{0}^{-1}z_{0}\delta
^{-1}z_{0}^{-1}g_{0}$.
By composition with the inner automorphism $x\rightarrow
z_{0}^{-1}g_{0}xg_{0}^{-1}z_{0},$ we get a homomorphism $\psi $ such
that $
\psi (\delta )=\delta ^{-1}$ for $\delta \in V$. The group generated by
$V$
is an open abelian subgroup of $\Gamma _{0}$.
\end{proof}

For any neighborhood $V$ of the identity,
$\gamma _{j}\subset V$ for large enough~$j$.
Hence there
exists an open abelian subgroup of $\Gamma _{0}$
if and only if $\gamma _{j}$ is abelian for some $j.$

\begin{corollary}\label{condition}
Let $M_{\infty }$ be a solenoid and let $K_\infty$ be the kernel
of $\pi_1(M_0)$.
Then $M_\infty$ is algebraically bihomogeneous if and only if
$\pi_1(M_j)/K_\infty$ is abelian for sufficiently large index~$j$, or,
equivalently, $\Gamma_j$ is abelian for sufficiently large index~$j$.
\end{corollary}

\section{An application}

Our algebraic condition for (topological)
bihomogeneity in Corollary \ref{condition}
is necessary but not sufficient.
For this, there should exist a homeomorphism $h\colon M_i\to M_i$
which induces an isomorphism $h_*\colon\pi_1(M_i)\to\pi_1(M_i)$
such that $h_*(x)=x^{-1}$ (modulo $K_\infty$).
The problem whether homomorphisms between fundamental groups
are realized by continuous maps is known as the geometric
realization problem.
It is a classical result of Nielsen \cite{Nielsen} that closed surfaces
admit geometric realizations.
This can be extended to certain three-dimensional
manifolds~\cite{Waldhausen}.
The following result now follows from Nielsen's theorem.

\begin{theorem}
A two-dimensional solenoid $S_\infty$
with kernel $K_\infty\subset \pi_1(S_0) $ is
bihomogeneous if and only if $\pi_1(S_i)/K_\infty$
is abelian for sufficiently large index $i$.
\end{theorem}

One easily constructs two-dimensional solenoids that are
not bihomogeneous, using results from geometric group theory.
The fundamental group $\pi_1(S)$ of a closed surface is subgroup separable,
see~\cite{Scott};
i.e., for every subgroup $H\subset \pi_1(S)$
there is a descending chain of subgroups of finite index with kernel $H$.
Hence, there
exists a solenoid with base-space $S$ and kernel $H$.
For a closed surface $S$ of genus greater than~$1$,
the fundamental group contains no abelian subgroup of finite index.
Therefore, a solenoid with base-space $S$
and kernel~$\{e\}$ is a (simply-connected) continuum which is
not bihomogeneous.

\section{Final remarks}

One-dimensional solenoids
are indecomposable continua.
It is not difficult to show that
higher-dimensional solenoids are not.
Rogers~\cite{Rogers2} has
shown
that a homogeneous, hereditarily indecomposable continuum is at most
one-dimensional. His question whether there exists a
homogeneous, indecomposable continuum of dimension greater than one
remains open.
\medbreak
Our example of a non-bihomogeneous space
is based on obstructions of the fundamental
group, which seems to be
characteristic for all examples so far.
So it is natural
to ask whether there exists a
simply-connected Peano
continuum that is homogeneous but not bihomogeneous. More generally, 
it is natural to ask whether there exists a continuum with trivial first \v{C}ech cohomology
that is homogeneous but not bihomogeneous.

\Addresses\recd

\end{document}